\documentclass[leqno,12pt]{amsart}
\usepackage{latexsym}
\usepackage{mathrsfs}
\usepackage{amsmath}
\usepackage{amssymb}
\usepackage[bookmarks]{hyperref}
\usepackage{pstricks,pst-plot}

\newtheorem{theorem}{Theorem}[section]
\newtheorem{lemma}[theorem]{Lemma}
\newtheorem{corollary}[theorem]{Corollary}
\newtheorem{proposition}[theorem]{Proposition}

\theoremstyle{definition}
\newtheorem{remark}[theorem]{Remark}
\newtheorem{definition}[theorem]{Definition}

\newcommand{\C}{\mathbb{C}}

\newcommand{\Z}{\mathbb{Z}}

\newcommand{\R}{\mathbb{R}}

\newcommand{\bthm}{\begin{theorem}}
\newcommand{\ethm}{\end{theorem}}
\newcommand{\blem}{\begin{lemma}}
\newcommand{\elem}{\end{lemma}}
\newcommand{\bcor}{\begin{corollary}}
\newcommand{\ecor}{\end{corollary}}
\newcommand{\bprop}{\begin{proposition}}
\newcommand{\eprop}{\end{proposition}}
\newcommand{\bdefn}{\begin{definition}}
\newcommand{\edefn}{\end{definition}}
\newcommand{\bpf}{\begin{proof}}
\newcommand{\epf}{\end{proof}}

\font\tenscr=rsfs10  scaled \magstep1

\font\sevenscr=rsfs7  scaled \magstep1
\font\fivescr=rsfs5  scaled \magstep1
\skewchar\tenscr='177 \skewchar\sevenscr='177 \skewchar\fivescr='177
\newfam\scrfam \textfont\scrfam=\tenscr \scriptfont\scrfam=\sevenscr
\scriptscriptfont\scrfam=\fivescr
\def\scr{\fam\scrfam}

\def \sm {\setminus}

\def\vp{\varphi}

\def\h#1{\widehat {#1}}

\def\ol {\overline}

\def\bC{\mathbb C}
\def\cn{\bC^n}
\def\YG{Y\cup \Gamma}
\def\what{\widehat}

\def\vp{\varphi}

\def\Rm{\begingroup\let\par=\null\obeylines\RemovE}
\def\RemovE#1\mR{\endgroup}

\def\sC{{\scr C}}

\begin{document}
\title[Arcs in Polynomially Convex Simple Closed Curves]{Polynomially Convex  Arcs in\\ Polynomially Convex Simple Closed Curves}
\author{Alexander J. Izzo}
\address{Department of Mathematics and Statistics, Bowling Green State University, Bowling Green, OH 43403}
\email{aizzo@bgsu.edu}
\thanks{The first author was supported by NSF Grant DMS-1856010.}
\author{Edgar Lee Stout}
\address{Department of Mathematics, University of Washington, Seattle, WA  98195}
\email {Stout@math.washington.edu}


\begin{abstract}
We prove that every polynomially convex arc is contained in a polynomially convex simple closed curve.  We also establish results about polynomial hulls of arcs and curves that are locally rectifiable outside a polynomially convex subset.
\end{abstract}
\maketitle

%
%

\section{The Results}

The main purpose of this paper is to prove the following theorem.  

\begin{theorem}
\label{1.23.20.i}
A polynomially convex arc $\lambda$ in $\cn$, $n\geq 2$, is contained in a polynomially convex simple closed curve $\gamma$ that can be chosen to lie in an arbitrarily small neighborhood of the given arc.  Furthermore, $\gamma$ can be chosen such that the open arc $\gamma\sm\lambda$ is $\sC^\infty$-smooth.  With this choice, if ${\mathcal P}(\lambda)={\mathcal C}(\lambda)$, then ${\mathcal P}(\gamma)={\mathcal C}(\gamma)$.
\end{theorem}

Recall that an \emph{arc} is a space homeomorphic to the closed unit interval, and a \emph{simple closed curve} is a space homeomorphic to the unit circle. 
As usual 
we denote by
${\mathcal C}(X)$ the space of all continuous $\bC$--valued functions on $X$ and by 
${\mathcal P}(X)$ the subspace of ${\mathcal C}(X)$ comprising those functions that can be approximated uniformly on $X$ by holomorphic polynomials.  Throughout the paper neighborhoods will be assumed to be open except when explicitly taken to be compact instead.

It is perhaps worth noting here that it remains an open question whether for every polynomially convex arc $\lambda$ in $\cn$,  the algebras $\mathcal C(\lambda)$ and $\mathcal P(\lambda)$ coincide.

The special case of the above theorem in which the arc is rectifiable was proved by the authors earlier, and in that case the simple closed curve can be chosen to be rectifiable as well \cite[Corollary~3.6]{Izzo-Stout:2022}.  It was also shown that when the arc is smooth of class $\sC^s$, $1\leq s \leq \infty$, then the closed curve can be chosen to be of class $\sC^s$ as well \cite[Corollary~3.9]{Izzo-Stout:2022}.  However, the methods used in the rectifiable and smooth cases are not applicable to nonrectifiable arcs, and the proof we will give of Theorem~\ref{1.23.20.i} is quite different from the proofs in 
\cite{Izzo-Stout:2022}.

As a first step toward proving Theorem~\ref{1.23.20.i}, we will prove the following  purely topological fact, which we suppose to be known but for which we have no reference. 

\begin{theorem}\label{12.27.20.i}
If $\lambda$ is an arc in $\R^n$, $n\geq 2$, and if $\Omega$ is an open set of 
$\R^n$ containing $\lambda$, then there is a simple closed curve $\gamma$ contained in $\Omega$ that contains $\lambda$ and is such that the open arc 
$\gamma\setminus\lambda$ is $\sC^\infty$-smooth.
\end{theorem} 

The case in which $\lambda$ is assumed to be rectifiable or smooth and a rectifiable or smooth simple closed curve containing it is sought was established by the authors in \cite{Izzo-Stout:2022}.
It should be noted that in the general situation of Theorem~\ref{12.27.20.i}, although the open arc $\gamma\setminus\lambda$ is smooth, the construction we will give does not guarantee that this open arc has finite length. 

As another step toward proving Theorem~\ref{1.23.20.i} we will establish the following polynomial convexity result.  
Here we use the standard notation that the polynomial hull of a compact set $X$ is denoted by $\what X$.

\begin{theorem}\label{arcs-in-lines}
Let $Y\subset \bC^n$ be a compact polynomially convex set, and let $\Gamma$ be an arc both of whose end points lie in $Y$ but that is otherwise disjoint from $Y$. Assume in addition that $\Gamma\setminus Y$ is locally rectifiable. Let $\Omega$ be a Stein neighborhood of $\widehat{Y\cup \Gamma}$, and let $V$ be a purely one-dimensional analytic subvariety of $\Omega$. If a nonempty open subarc of $\Gamma$ is contained in $V$, then either $Y\cup \Gamma$ is polynomially convex or else $\what{(Y\cup \Gamma)}\setminus Y$ is contained in $V$.  The same conclusion holds with the arc $\Gamma$ replaced by a simple closed curve that intersects $Y$ in a single point.
\end{theorem}

The following corollary is immediate.

\bcor\label{needed-case}
Under the hypotheses of Theorem~\ref{arcs-in-lines}, if a nonempty open subarc of 
$\Gamma$ is contained in $V$ but $\Gamma$ is not entirely contained in $V$, then $Y\cup \Gamma$ is polynomially convex.
\ecor

For instance if a nonempty open subarc of a nonpolynomially convex rectifiable simple closed curve is contained in a purely one-dimensional analytic variety, then the entire simple closed curve must be contained in that variety.

We will establish the case of Theorem~\ref{arcs-in-lines} in which $\Gamma$ is rectifiable using  a standard result \cite[Theorem~3.1.1]{Stout2007} about polynomial convexity and sets of finite length which we state here for the reader's 
convenience\footnote{As stated in \cite{Stout2007}, the hypothesis in the last sentence of the theorem is that the map $\check H^1(\YG; \Z)\rightarrow \check H^1(Y;\Z)$ is an isomorphism.  However, as is evident from the proof given there, it is enough to assume that this map is a monomorphism.  Also the assertion in the last sentence that $Y\cup \Gamma$ is polynomially convex is not stated in the theorem in \cite{Stout2007} but is established in the proof given there.}.

\begin{theorem}\label{311}
Let $Y$ be a compact polynomially convex subset of $\cn$, and let $\Gamma$ be a subset of $\cn$ contained in a compact connected set of finite length such that $Y\cup\Gamma$ is compact.  The polynomial hull $\what{Y\cup\Gamma}$ has the property that the complementary set $\what{Y\cup\Gamma}\setminus (Y\cup\Gamma)$ either is empty or else is a purely one-dimensional analytic subvariety of $\cn\setminus(Y\cup\Gamma)$. If the map $\check H^1(\YG; \Z)\rightarrow \check H^1(Y;\Z)$ induced by the inclusion $Y\hookrightarrow \YG$ is a monomorphism, then
$Y\cup\Gamma$ is polynomially convex and ${\mathcal P}(\YG)=\{ f\in {\mathcal C}(\YG): f|Y\in {\mathcal P}(Y)\}$. 
\end{theorem}

To obtain the general case of Theorem~\ref{arcs-in-lines} in which $\Gamma\setminus Y$ is only \emph{locally} rectifiable, we will prove the following partial generalization of the above theorem.

\begin{theorem}\label{general311}
Let $Y$ be a compact polynomially convex subset of $\cn$, and let $\Gamma$ be a subset of $\cn$ such that $\YG$ is compact and such that for every neighborhood $U$ of $Y$ in $\cn$, the set $\Gamma\setminus U$ is contained in a compact connected set of finite length.  Suppose also that the map $\check H^1(\YG; \Z)\rightarrow \check H^1(Y;\Z)$ induced by the inclusion $Y\hookrightarrow \YG$ is a monomorphism.  Then $\YG$ is  polynomially convex and ${\mathcal P}(\YG)=\{ f\in {\mathcal C}(\YG): f|Y\in {\mathcal P}(Y)\}$.
\end{theorem}

It is well known (and immediate from Theorem~\ref{311}) that every rectifiable arc in $\cn$ is polynomially convex and that for a rectifiable simple closed curve $\gamma$ in $\cn$, the set $\what\gamma \setminus \gamma$, if not empty, is a purely one-dimensional analytic subvariety of 
$\cn\setminus\gamma$.  As applications of Theorem~\ref{general311} we obtain generalizations in which the arc or simple closed curve is required to be only \emph{locally} rectifiable outside a polynomially convex subset.

\begin{theorem}\label{locally-rectifiable-arc}
Let $\lambda$ be an arc in $\cn$, and let $E$ be a compact subset of $\lambda$ that is polynomially convex.  Suppose that $\lambda\setminus E$ is locally rectifiable.  Then $\lambda$ is polynomially convex.  Furthermore, if ${\mathcal P}(E)={\mathcal C}(E)$, then ${\mathcal P}(\lambda)={\mathcal C}(\lambda)$.
\end{theorem}

As special cases we mention the following immediate corollaries.
(Note that every compact set $E$ of zero length is polynomially convex and satisfies 
${\mathcal P}(E)={\mathcal C}(E)$ \cite[Theorem~1.6.2]{Stout2007}.)

\begin{corollary}\label{cor1}
If an arc $\lambda$ in $\cn$ is locally rectifiable off a closed subset of  zero length, then $\lambda$ is polynomially convex and satisfies ${\mathcal P}(\lambda)={\mathcal C(\lambda)}$.
\end{corollary}

\begin{corollary}\label{cor2}
An arc $\lambda$ in $\cn$ whose interior is locally rectifiable is polynomially convex and satisfies $\mathcal P(\lambda)=\mathcal C(\lambda)$.
\end{corollary}

By analogy with Theorem~\ref{311},
one might be tempted to imagine that without the 
hypothesis regarding cohomology in Theorem~\ref{general311}, the set $\widehat {(\YG)}\setminus(\YG)$, if nonempty, must be a purely one-dimensional analytic subvariety of $\cn\setminus(\YG)$.  However, this is not so, as is shown by an example of Herbert Alexander
\cite{Alexander1986} of a set that is a countable union of rectifiable simple closed curves whose polynomial hull is not an analytic variety but is instead a countable union of analytic varieties.  Alexander's example is a union of the circle $C=\{(z,0)\in \bC^2: |z|=1\}$ and a countable collection $\{\gamma_k\}_{k=1,2,\ldots}$ of analytic simple closed curves that accumulate on $C$.  Thus setting $Y=\what C$ and $\Gamma=\cup_{k=1,2,\ldots} \gamma_k$, for every neighborhood $U$ of $Y$ the set $\Gamma\setminus U$ is contained in a  finite union of simple closed curves of finite length and consequently is contained in a connected set of finite length.

The proof of Theorem~\ref{general311} does, however, show that in the setting of the theorem with the cohomology hypothesis removed, the set $Y\cup \Gamma$ satisfies Gabriel Stolzenberg's generalized argument principle (which we recall in Section~\ref{proofs-section}).  Furthermore from Theorem~\ref{general311} we will obtain the following analogue of Theorem~\ref{locally-rectifiable-arc} for simple closed curves, thus generalizing the well-known result about polynomial hulls of rectifiable simple closed curves.

\begin{theorem}\label{closed-curve-cor}
Let $\gamma$ be a simple closed curve in $\cn$, and let $E$ be a compact subset of $\gamma$ that is polynomially convex.  Suppose that $\gamma\setminus E$ is locally rectifiable.  Then ${\what\gamma} \setminus \gamma$ either is empty or else  is a purely 
one-dimensional analytic subvariety of $\cn\setminus\gamma$.
\end{theorem}

As a corollary we obtain the following.

\begin{corollary}\label{cor3}
If a simple closed curve $\gamma$ in $\cn$ is locally rectifiable off a closed subset of zero length, then ${\what\gamma} \setminus \gamma$ either is empty or else  is a purely 
one-dimensional analytic subvariety of $\cn\setminus\gamma$.  In the former case, ${\mathcal P}(\gamma)={\mathcal C(\gamma)}$.
\end{corollary}

We prove the topological Theorem~\ref{12.27.20.i} in the next section.  The proofs of our other results are given in Section~3, except for the proof of Theorem~\ref{1.23.20.i} which we give here.
\medskip

\noindent{\bf Proof of Theorem \ref{1.23.20.i}.} Let $U$ be an arbitrary neighborhood of the polynomially convex arc $\lambda$.  By Theorem~\ref{12.27.20.i} there is a simple closed curve $\gamma$ contained in $U$ that contains $\lambda$ and is such that $\gamma\setminus\lambda$ is $\sC^\infty$-smooth.
By Corollary~\ref{needed-case} (with $Y=\lambda$, $\Gamma$ the closure of the open arc $\gamma\setminus\lambda$, $\Omega=\cn$, and $V$ a complex line), all that is needed to see that $\gamma$ can be chosen so as to be polynomially convex is to note that we can modify $\gamma$ so that some open subarc of $\gamma\setminus\lambda$ is contained in a complex line but $\gamma\setminus\lambda$ is not entirely contained in the complex line.
That ${\mathcal P}(\lambda)={\mathcal C}(\lambda)$ then implies ${\mathcal P}(\gamma)={\mathcal C}(\gamma)$ follows from the following well-known result
\cite[Corollary~1.6.8]{Stout2007}.

\begin{theorem}\label{168}
If $Y$ is a rationally convex subset of $\cn$ and $\Gamma$ is a subset of $\cn$ of zero two-dimensional Hausdorff measure such that $K=Y\cup\Gamma$ is compact, then $K$ is rationally convex, and ${\mathcal R}(K)=\{f\in {\mathcal C}(K): f|Y\in {\mathcal R}(Y)\}$.
\end{theorem}

Again our notation is standard: $\mathcal R(K)$ is the algebra of all functions that can be approximated uniformly on $K$ by restrictions to $K$ of rational functions holomorphic on a neighborhood of $K$.

\begin{remark}
In fact the construction in the proof of Theorem~\ref{12.27.20.i} can easily be carried out so as to yield a polynomially convex simple closed curve $\gamma$ at the outset so that no modification of the curve $\gamma$ is needed in the proof of Theorem~\ref{1.23.20.i}.  To see this, simply observe that the curve $\gamma$ that we construct contains 
straight line segments and that obviously the points $b_k$, $k=0,1,2,\ldots$ selected in the construction, and hence the open arc $\gamma\setminus\lambda$, can be chosen so as to be contained in no proper vector subspace of $\R^n$.
\end{remark}


%
%

\section{Proof of the Topological Theorem~\ref{12.27.20.i}}

We first treat the case that $n\geq 3$; the planar case, which is not needed for the proof of Theorem \ref{1.23.20.i}, requires a different argument.
Throughout, by \emph{smooth} we mean of class $\sC^\infty$.

 With no loss of generality, suppose $\Omega$ to be connected. Let the end points of $\lambda$ be the points $a$ and $a'$. Let $B_0$ be an open  ball contained in $\Omega$ and centered at $a$, and let $B_0'$ be an open ball contained in $\Omega$ and centered at $a'$. 
Choose the balls $B_0$ and $B_0'$ small enough that their closures are disjoint and are contained in $\Omega$.  For notational convenience we suppose the balls $B_0$ and $B_0'$ each to be of radius one.  
For $k=1,2,\dots$, let $B_k$ be the open ball of radius $2^{-k}$ centered at the point $a$, and let $B_k'$ be the corresponding ball centered at the point $a'$. For each $k=0,1,2,\ldots$, let $b_k$ be a point of $bB_k\setminus \lambda$ and $b_k'$ a point of $bB_k'\setminus\lambda$.  (Here $bB_k$ denotes the boundary of $B_k$.)

The set $\Omega\setminus(B_0\cup B_0')$ is a connected manifold-with-boundary.  Since a space of topological dimension one cannot disconnect a manifold-with-boundary of dimension greater than or equal to three (see \cite[p.~48]{Hurewicz-Wallman:1948}), the set $\Omega\setminus(B_0\cup B_0'\cup\lambda)$ is also a connected manifold-with-boundary.
Therefore, there is a smooth arc $\rho$ in $\Omega\setminus(B_0\cup B_0'\cup\lambda)$ from the point $b_0$ of $bB_0$ to the point $b_0'$ of $bB_0'$.
In addition, we can choose $\rho$ so that its interior is disjoint from $bB_0\cup bB_0'$ and so that a subarc of $\rho$ with end point $b_0$ lies in the real line through the points $a$ and $b_0$ and a subarc of $\rho$ with end point $b_0'$  lies in the real line through the points $a'$ and $b_0'$.

For each $k=0,1,2,\ldots$, the set $\ol B_k\setminus(B_{k+1}\cup \lambda)$ is a connected 
manifold-with-boundary (again by \cite[p.~48]{Hurewicz-Wallman:1948}), so there is a smooth arc $\tau_k$ in this set from $b_k$ to $b_{k+1}$.  In addition, we can choose $\tau_k$ so that its interior is disjoint from $bB_k\cup bB_{k+1}$ and so that a subarc of $\tau_k$ with end point $b_k$ lies in the real line
through the points $a$ and $b_k$ and a subarc of $\tau_k$ with end point $b_{k+1}$ lies in the real line through the points $a$ and $b_{k+1}$.
Then the union $\tau=(\cup_{k=0,1,2,\dots}\tau_k)\cup\{a\}$ is an arc  from $b_0$ to $a$ whose interior lies in $B_0$.  Also $\tau$ intersects $\lambda$ only in the point $a$. The 
arc $\tau$ is smooth except possibly at the end point $a$. 

By a similar construction we obtain an arc $\tau'$ from $b_0'$  to $a'$ such
that the interior of $\tau'$ lies in $B_0'$, such that $\tau'$ intersects $\lambda$ only in the point $a'$, and such that $\tau'$ is smooth except possibly at the end point $a'$. 

The union $\gamma=\lambda\cup \rho\cup\tau\cup\tau'$
is a simple closed curve in $\Omega$ that contains $\lambda$.
Furthermore, the open arc $\gamma\setminus\lambda$ is smooth.

The theorem is proved in the case of arcs in $\R^n$ with $n$ at least three. The argument above does not work in the case of arcs in the plane.  For that case we will apply conformal mapping methods.

The following well known theorem is proved in the paper of Philip Church \cite[Theorem~2.1]{
Church:1964}.
In fact, the theorem seems to be due to Marie Torhorst \cite{Torhorst:1917}. 

 \begin{theorem} If $G$ is a simply connected region on the Riemann sphere whose boundary is a nondegenerate Peano continuum, then each prime end of $G$ is a single point. Thus,
if $\vp$ is a conformal homeomorphism of the open unit disc  $\mathbb{U}$ onto $G$, then $\vp$ extends to
$\overline{\mathbb {U}}$ as a continuous map.
\end{theorem}

With $\bC^*$ the Riemann sphere, there is a conformal map $\vp:\mathbb {U}
\rightarrow\bC^*\setminus \lambda$ of $\mathbb {U}$ onto $\bC^*\setminus \lambda$.  By the Church-Torhorst theorem above, $\vp$ extends continuously to the closure of $\mathbb {U}$.
With  $a$ and $a'$ the end points of $\lambda$, choose points $p$ and $p'$ in the boundary of $\mathbb {U}$ that satisfy $\vp(p)=a$ and $\vp(p')=a'$. For $\ell$ an arc in $\mathbb{U}\cup\{p, p'\}$ with endpoints $p$ and $p'$ and interior contained in $\mathbb {U}$, setting 
$\gamma=\lambda\cup\vp(\ell)$ we obtain a simple closed curve in the plane that contains $\lambda$.
By choosing $\ell$ to be smooth and to lie sufficiently near $b\mathbb {U}$, we can insure that the open arc $\gamma\setminus\lambda$ is smooth and contained in $\Omega$. 

\medskip

\begin{remark}
In fact, the smoothness conclusion in Theorem~\ref{12.27.20.i} can be strengthened: the simple closed curve $\gamma$ can be chosen so that the open arc $\gamma\setminus\lambda$ is real-analytic.  To see this, choose a mapping $G:\R\rightarrow\R^n$ that takes the real line diffeomorphically onto the open arc $\gamma\setminus\lambda$ and satisfies $\lim_{x\rightarrow-\infty}G(x)=a$ and $\lim_{x\rightarrow\infty}G(x)=a'$.  By the openness of $\sC^1$ embeddings in the strong topology (also know as the fine or Whitney topology) on 
$\sC^1(\R,\R^n)$ \cite[Chapter~2, Theorem~1.4]{Hirsch:1976}, there is a positive continuous function $\delta$ on $\R$ such that every $\R^n$-valued $\sC^1$ mapping $F$ on $\R$ that satisfies
\begin{equation}\label{embedding-condition}
|F(x)-G(x)|< \delta(x)\ \ {\rm and\ \ } |F'(x)-G'(x)|<\delta(x) \ \ {\rm for\ all\ } x\in \R
\end{equation}
is an embedding.
Since we can replace $\delta$ by a smaller positive continuous function, we may assume that $\delta$ tends to zero at $\pm \infty$ and that for every $x\in \R$, the number $\delta(x)$ is smaller than the distance from $G(x)$ to the closed set $\C^n\setminus (\Omega\setminus\lambda)$.  Then every mapping $F$ satisfying (\ref{embedding-condition}) has range in $\Omega\setminus\lambda$ and satisfies $\lim_{x\rightarrow-\infty}F(x)=a$ and $\lim_{x\rightarrow\infty}F(x)=a'$.  Thus to establish the existence of the desired simple closed curve $\gamma$ with $\gamma\setminus\lambda$ real-analytic, it suffices to show that there exists a real-analytic mapping $F$ that satisfies (\ref{embedding-condition}).

Recall that $\R$ is a {\it{Carleman continuum}}, i.e.,  that given a continuous $\C$--valued function $g$ on $\R$, and regarding $\R$ as the real axis in the complex plane $\C$, for every positive continuous function $\varepsilon$ on $\R$ there is an entire function $f$ on $\C$ such that 
$$|f(x)-g(x)|<\varepsilon(x)\ \ {\rm for\ all\ } x\in \R.$$
Indeed, more is true: If $g$ is of class $\sC^k$, $k$ a positive integer, then there is an entire function $f$ on $\C$ with 
\begin{equation}\label{ck-approx}
|f^{(j)}(x)-g^{(j)}(x)|<\varepsilon(x)\ \ \hbox{for all $x\in\R$ and all $j=0,1,\dots,k$}.
\end{equation}
The existence of the function $f$ was established in the case of first derivatives by Wilfred Kaplan \cite[Theorem~3]{Kaplan:1955}.  The result for derivatives of higher order (which we do not need) was given by Lothar Hoischen \cite[Satz~2]{Hoischen:1973}.  Of course when $g$ is $\R$--valued, (\ref{ck-approx}) continues to hold with $f$ replaced by the real part of $f$.  Applying this approximation result (with $k=1$) to each component of $G$ yields the desired real-analytic mapping $F$.
\end{remark}

\section{Proofs of the Polynomial Convexity Results}\label{proofs-section}

Our results regarding polynomial convexity (with the exception of Theorem~\ref{1.23.20.i} which was proved in the introduction) will be proved in the following order: Theorem~\ref{general311}, Theorem~\ref{arcs-in-lines}, Theorem~\ref{locally-rectifiable-arc}, Theorem~\ref{closed-curve-cor}, Corollary~\ref{cor3}.

To prove Theorem~\ref{general311} we will use Stolzenberg's generalized argument principle \cite{Stol1}.

\bdefn
A compact set $X$ in $\cn$ satisfies the {\em generalized argument principle\/} provided that, if $P$ is a polynomial that has a continuous logarithm on $X$, then $0\notin P(\h X)$.
\edefn

Theorem~\ref{general311} will be proved by showing that $\YG$ satisfies the generalized argument principle and then obtaining the result from Theorem~\ref{311} 
and the construction of a certain polynomial that has a continuous logarithm on $\YG$.  Note that we actually need Theorem~\ref{general311} only in the case when $Y$ and $\YG$ are arcs.  In that case the proof simplifies; once we have shown that $\YG$ satisfies the generalized argument principle we just need to note that $\YG$ is rationally convex by Theorem~\ref{168}  and then invoke the following observation of Stolzenberg.

\bprop \cite[Corollary 2.3]{Stol1}
If $X$ is simply coconnected, i.e., satisfies $\check H^1(X;\mathbb Z)=0$, and satisfies the generalized argument principle, then the polynomial and rational hulls 
of $X$ coincide.
\eprop

The following modification of \cite[(2.4)]{Stol1} is stated in \cite[Proposition~3.5]{Izzo2019}; it follows from the classical argument principle in the same manner as \cite[(2.4)]{Stol1}.

\bprop \label{classicalarg}
If $V$ is a purely one-dimensional analytic subvariety of an open subset of $\cn$ with $\ol V$ compact, and $P$ is a polynomial with a continuous logarithm on $bV=\ol V\setminus V$, then $0\notin P(\ol V)$.
\eprop

Thus the boundary of a relatively compact purely one-dimensional analytic variety $V$ with $\h {bV}=\ol V$ satisfies the generalized argument principle.

The following result of Stolzenberg shows that the generalized argument principle is preserved by certain limits.

\bthm \cite[(2.12)]{Stol1}  \label{limit}
Let $(X_k)_{k=1,2,\dots}$ be a sequence of compact sets in $\cn$ each of which satisfies the generalized argument principle.  If in the Hausdorff metric $X_k\rightarrow X$ and $\h X_k\rightarrow \h X$, then $X$ satisfies the generalized argument principle.
\ethm

\blem\label{Hausdorffconvergence}
If $(X_k)_{k=1,2,\dots}$ is a decreasing sequence of compact sets in $\cn$ such that $\bigcap_{k=1,2,\ldots} X_k=X$, then in the Hausdorff metric $X_k\rightarrow X$ and $\h X_k\rightarrow \h X$.
\elem

\noindent{\bf{Proof.}}
It follows from the hypotheses of the lemma that $\bigcap_{k=1,2,\ldots}\widehat X_k= \widehat X$.  
Since  for every decreasing sequence $(Y_k)_{k=1,2,\dots}$ of compact sets in $\cn$ one has that $Y_k\rightarrow \bigcap_{k=1,2,\dots} Y_k$, the lemma follows.

\medskip

\noindent{\bf{Proof of Theorem~\ref{general311}.}}
Set $X=\YG$.  We first show that $X$ satisfies the generalized argument principle.
Choose a decreasing sequence of compact polynomially convex neighborhoods 
$(L_k)_{k=1,2,\dots}$ of $Y$ with intersection $Y$.  Let $X_k=L_k\cup \Gamma$.  Then $\bigcap_{k=1,2,\ldots} X_k=X$.  By hypothesis $\Gamma\setminus L_k$
is contained in a compact connected set of finite length.  Thus by Theorem~\ref{311}, $\what X_k\setminus X_k$ either is empty or else is a purely one-dimensional analytic subvariety of $\cn\setminus X_k$.

If $X_k$ is polynomially convex, then it satisfies the generalized argument principle. If $\what X_k\setminus X_k$ is a variety, say $V$, then a polynomial  that has a logarithm on $X_k$ has a logarithm on $bV$ and so, by Proposition~\ref{classicalarg}, has no zero on $V$. Thus, $X_k$ is again seen to satisfy the generalized argument principle. 
  By Lemma~\ref{Hausdorffconvergence}, in the Hausdorff metric $X_k\rightarrow X$ and $\what X_k\rightarrow \what X$.  Thus by Theorem~\ref{limit}, $X$ satisfies the generalized argument principle.

To conclude the proof that $X$ is polynomially convex, we assume that $\what X\setminus X$ is nonempty and derive a contradiction to the generalized argument principle.  (For this we essentially follow the beginning of the proof of Theorem~\ref{311} given in \cite[p.~151]{Stout2007}.)  We can assume without loss of generality that $\Gamma$ is disjoint from $Y$.  Then $\Gamma$ is a countable union of sets of finite length and hence has two-dimensional Hausdorff measure zero.  Consequently, $X$ is rationally convex by 
Theorem~\ref{168}.  Now fix a point $p\in \what X \setminus X$.  By the polynomial convexity of $Y$, there exists a polynomial $Q$ such that $\Re Q<0$ on $Y$ and $Q(p)=0$.  By the rational convexity of $X$, there exists a polynomial $R$ such that $0\notin R(X)$ and $R(p)=0$.  Choose a positive constant $c$ sufficiently large that $\Re(cQ+ \zeta R)<0$ on $Y$ for all choices of $\zeta\in \bC$ with 
$|\zeta|<1$.  No matter what the choice of $c$ and $\zeta$, we have $(cQ + \zeta R)(p)=0$. The function $H=cQ/R$ is holomorphic on a neighborhood of $X$, so the set $H(\Gamma)$ has zero area in $\bC$. Choose $\zeta$ such that 
$|\zeta|<1$ and $-\zeta\notin H(\Gamma)$.  Then the polynomial $P=cQ+\zeta R$ vanishes at $p$ but is zero-free on $X$.  Because $\Re P<0$ on $Y$, the function $P$ has a continuous logarithm on $Y$.  The hypothesis that the map $\check H^1(X; \Z)\rightarrow \check H^1(Y;\Z)$ induced by inclusion
is a monomorphism then implies that $P$ has a continuous logarithm on $X$.  We thus have a contradiction to the generalized argument principle for $X$.

The final assertion of the theorem now follows from Theorem~\ref{168}.

\medskip

\noindent{\bf {Proof of Theorem~\ref{arcs-in-lines}.}}
We treat first the case that $\Gamma$ is rectifiable.

Suppose that some nonempty open subarc $\sigma$ of $\Gamma$ is contained in $V$ and that $Y\cup\Gamma$ is not polynomially convex.  We first show that $\what{Y\cup \Gamma}\setminus(Y\cup \Gamma)$ is contained in $V$ and then obtain that, in fact, $\what {Y\cup\Gamma}\setminus Y$ is contained in $V$.   

Fix a point $p\in\what {Y\cup \Gamma}\setminus(Y\cup\Gamma)$, and choose a Jensen measure $\mu$ for the functional of evaluation at $p$ on ${\mathcal P}(Y\cup\Gamma)$ supported on $Y\cup\Gamma$.
Note that $Y$ is a deformation retract of $Y\cup (\Gamma\setminus\sigma)$, and hence, the map $\check H^1\bigl(Y\cup(\Gamma\setminus\sigma); \Z\bigr)\rightarrow \check H^1(Y;\Z)$ induced by inclusion is an isomorphism.  Thus
$Y\cup(\Gamma\setminus\sigma)$ is polynomially convex by the last part of Theorem~\ref{311}, so necessarily $\mu(\sigma)>0$.  Therefore, for every holomorphic function $f$ on $\Omega$ that vanishes identically on $V$,
$$\log|f(p)|\leq \int\, \log|f|\, d\mu= -\infty,$$
so $f(p)=0$. Consequently, $p$ is in $V$, because since $V$ is an analytic variety in the Stein domain $\Omega$, for each point $q$ not in $V$ there is a holomorphic function 
$f$ on $\Omega$ that vanishes identically on $V$ but not at $q$ \cite[Theorem~I5]{Gunning:1990}.

We have shown that $\what{Y\cup \Gamma}\setminus( Y\cup\Gamma)$ is contained in $V$.  Let $W=\what{Y\cup\Gamma}\setminus( Y\cup\Gamma)$.  To conclude the proof it suffices to show that $\Gamma$ is contained in the closure of $W$.  By Theorem~\ref{311}, $W$ is an analytic subvariety of $\bC^n\setminus (Y\cup\Gamma)$.  The boundary $bW$ of $W$ lies in $Y\cup\Gamma$, and by the maximum principle, $W$ is contained in the polynomial hull of $bW$.  Since removing an arbitrary open subarc of $\Gamma$ from $Y\cup\Gamma$ yields a polynomially convex set (by the last part of Theorem~\ref{311}), we conclude that $bW$ must contain $\Gamma$ thereby concluding the proof in the special case that $\Gamma$ is rectifiable.

With the case that $\Gamma$ is rectifiable in hand, we turn to the general case, i.e., the case that $\Gamma\setminus Y$ is locally rectifiable, and hence, every compact subset of $\Gamma\setminus Y$ has finite length.
Theorem~\ref{general311} shows that for $\sigma$ an open subarc of $\Gamma\setminus Y$ whose closure is contained in $\Gamma\setminus Y$, the compact set $X=Y\cup(\Gamma\setminus\sigma)$ is polynomially convex.  Note that $X\cup\sigma=Y\cup\Gamma$.  The closure of $\sigma$ is a rectifiable arc with end points in $X$ and otherwise disjoint from $X$.  Thus the special case of the theorem already established (applied with $X$ in place of $Y$) shows that either $X\cup\sigma=Y\cup\Gamma$ is polynomially convex or else $\what{X\cup\sigma}\setminus X=\what{Y\cup\Gamma}\setminus X$ is contained in $V$.  Since each point of $\Gamma$ is contained in 
$\what{Y\cup\Gamma}\setminus X$ for some choice of $\sigma$, this shows that
$\what{Y\cup\Gamma}\setminus Y$ is contained in $V$ if $Y\cup\Gamma$ is not polynomially convex.
\medskip

\noindent{\bf Proof of Theorem \ref{locally-rectifiable-arc}.} 
This follows immediately from Theorem~\ref{general311} by setting $Y=E$ and $\Gamma=\lambda$, and noting that for every neighborhood $U$ of $Y$, the set 
$\Gamma\setminus U$ is contained in the union of finitely many rectifiable arcs, e.g., the union of finitely many components of $\Gamma\setminus U'$ for $U'$ a neighborhood of $Y$ whose closure is contained in $U$.

\medskip

\noindent{\bf Proof of Theorem \ref{closed-curve-cor}.} 
Let $\sigma$ be an open subarc of $\gamma\setminus E$ whose closure $\Gamma$ is an arc of finite length contained in $\gamma\setminus E$.  Then the arc $Y=\gamma\setminus \sigma$ is polynomially convex by Theorem~\ref{general311}.  The result now follows immediately by applying Theorem~\ref{311} to $\gamma=Y\cup\Gamma$.

\medskip

\noindent{\bf Proof of Corollary \ref{cor3}.} 
The first assertion is immediate from Theorem~\ref{closed-curve-cor}.  The second assertion follows from Theorem~\ref{168}.

\end{document}